\documentclass{amsart}
\usepackage{epsfig}
\usepackage{amsmath}
\usepackage{amsfonts}
\begin{document}

\newtheorem{thm}{Theorem}
\newtheorem{lem}[thm]{Lemma}
\newtheorem{cor}[thm]{Corollary}
\newtheorem{conj}[thm]{Conjecture}
\newtheorem{qn}{Question}
\newtheorem{pro}{Proposition}[section]

 \theoremstyle{definition}
\newtheorem{defn}{Definition}[section]

\theoremstyle{remark}
\newtheorem{rmk}{Remark}[section]

\def\square{\hfill${\vcenter{\vbox{\hrule height.4pt \hbox{\vrule width.4pt
height7pt \kern7pt \vrule width.4pt} \hrule height.4pt}}}$}

\def\Z{\Bbb Z}
\def\R{\Bbb R}
\newenvironment{pf}{{\it Proof:}\quad}{\square \vskip 12pt}

\title{Groups generated by two  elliptic elements in $\mathbf{PU}(2,1)$}
\author{Baohua Xie, Yueping Jiang }
\address{College of Mathematics and Econometrics \\ Hunan University \\ Changsha, 410082, China}
\email{xiexbh@gmail.com, ypjiang731@163.com}
\maketitle

\begin{abstract}
Let $f$ and $g$ be two elliptic elements in $\mathbf{PU}(2,1)$ of
order $m$ and $n$ respectively, where $m\geq n>2$. We prove that if
the distance $\delta(f,g)$ between the complex lines or points fixed
by $f$ and $g$ is large than a certain number, then the group
$\langle f, g\rangle$ is discrete nonelementary and isomorphic to
the free product $\mathbf{Z}_{m}*\mathbf{Z}_{n}$.
\end{abstract}

\section{Introduction}
 A subgroup of Fuchsian groups or Kleinian groups generated by two
 elements was studied by many authors. An interesting question is to explore the conditions for two
 elements in Fuchsian groups or Kleinian groups to generate discrete free
 group. In \cite{k}, Knapp found necessary and sufficient conditions
 for two elliptic transformations to generate a discontinuous
 subgroup of $Lf(2, \mathbf{R})$, the group of linear fractional
 transformations. Lyndon and Ullman \cite{lu} gave conditions for two
 hyperbolic transformations whose fixed points separate each other
 to generate a discrete free group of rank 2. In general, Purzitsky \cite{pu} find
 necessary and sufficient conditions for the subgroups generated by any
 pair $A,B \in Lf(2, \mathbf{R})$ to be the discrete free product of
 the cyclic groups $\langle A\rangle$ and $\langle B \rangle$.

The following Theorem 1 is well known in real hyperbolic geometry.
    It is essentially contained in \cite{k}.
\begin{thm}
   Suppose that $f$ and $g$ are elliptic elements of $\mathbf{PSL}(2,\mathbf{R})$  of order $m$ and $n$.
  Let $\delta(f,g)$ be the distance between the fixed points of $f$ and
  $g$. If

  $$\cosh\delta(f,g)>\frac{\cos\frac{\pi}{m}\cos\frac{\pi}{n}+1}{\sin\frac{\pi}{m}\sin\frac{\pi}{n}},$$
  then $\langle f,g\rangle$  is  discrete and isomorphic to the free product $\langle f\rangle \ast \langle g\rangle$.

\end{thm}

 F. W Gehring, C. Maclachlan and G. J.
Martin proved a similar result in the case of Kleinian groups.
\begin{thm} \cite{gmm}
Suppose that $f$ and $g$ are elliptic elements of $\mathbf{
PSL}(2,\mathbf{C})$  of order $m$ and $n$.
  Let $\delta(f,g)$ be the distance between the axes of $f$ and
  $g$. If

  $$\cosh\delta(f,g)>\frac{\cos\frac{\pi}{m}\cos\frac{\pi}{n}+1}{\sin\frac{\pi}{m}\sin\frac{\pi}{n}},$$
  then $\langle f,g\rangle$  is  discrete and isomorphic to the free product $\langle f\rangle \ast \langle g\rangle$.

\end{thm}
 In this paper, The principal problem we wish to consider is that of giving condition in terms of transformations in
 complex hyperbolic $2$-space for the free product of a cyclic group of order $p$ and a cyclic group of
 order $q$. We prove the following
result.

\begin{thm}
Let $f, g\in \mathbf{ PU}(2,1) $ be elliptic elements
 of order $m$ and $n$.  Let $f$ and $g$ be in one of
the following cases

 {\rm (1)} $f$ and
  $g$ are reflections in  complex lines;

{\rm (2)} $f$ is reflection in a complex line and $g$ is reflection
in a point;

{\rm (3)} $f$ and
  $g$ are reflections in  points.\\
Let $\delta(f,g)$ be the distance between the complex lines or
points fixed by  $f$ and $g$. Then
$$\cosh\delta(f,g)>\frac{\cos\frac{\pi}{m}\cos\frac{\pi}{n}+1}{\sin\frac{\pi}{m}\sin\frac{\pi}{n}}$$
 will imply  that $\langle f,g\rangle$  is  discrete and isomorphic to the free product $\langle f\rangle \ast \langle g\rangle$.

\end{thm}

The pattern of this result is very similar to the analogous results
in real hyperbolic space. A possible application of Theorem 3 is in
the study of complex hyperbolic triangle groups, see for example A.
Pratoussevitch \cite{pr}, R. E. Schwartz \cite{sc}.

\section{Complex hyperbolic space}

 First, we recall some terminology. More details can be found in
 \cite{bm,go,jkp,k1,k2}.
Let $\mathbf{C}^{2,1}$ denote the complex vector space of dimension
3, equipped with a non-degenerate Hermitian form  of signature
(2,1). There are several such forms. We use the following form,
called the second Hermitian form
$$\langle \mathbf{z}, \mathbf{w}\rangle=\mathbf{w}^{*} J\mathbf{z}$$
where $\mathbf{z}, \mathbf{w}$ are column vectors in
$\mathbf{C}^{2,1}$, the Hermitian transpose is denote by $.*$ and
$J$ is the Hermitian matrix
                      $$J= \begin{bmatrix} 0& 0& 1\\
0& 1& 0\\1&0& 0\end{bmatrix}.$$
 Consider the following subsets of $\mathbf{C}^{2,1}$

$$V_{+}=\{\mathbf{v}\in \mathbf{C}^{2,1}|\langle \mathbf{v},\mathbf{v}\rangle>0\},$$
$$V_{-}=\{\mathbf{v}\in \mathbf{C}^{2,1}|\langle \mathbf{v},\mathbf{v}\rangle<0\},$$
$$ V_{0}=\{\mathbf{v}\in \mathbf{C}^{2,1}|\langle
\mathbf{v},\mathbf{v}\rangle=0\}.$$

Let $\mathbf{P}:\mathbf{C}^{2,1}-\{0\}\rightarrow \mathbf{CP}^{2,1}$
be the canonical projection onto complex projective space. Then
$\mathbf{H}_{\mathbf{C}}^2=\mathbf{P}(V_{-})$ associated with the
Bergman metric is complex hyperbolic space. The biholomorphic
isometry group of $\mathbf{H}_{\mathbf{C}}^2$ is $\mathbf{PU}(2,1)$
acting by linear projective transformations. Here $\mathbf{PU}(2,1)$
is the projective unitary group with respect to the Hermitian form
defining on $\mathbf{C}^{2,1}$. In other words, for all $\mathbf{z}$
and $\mathbf{w}$ in $\mathbf{C}^{2,1}$ we have

$$\mathbf{w}^{*} J\mathbf{z}=\langle \mathbf{z}, \mathbf{w}\rangle=\langle B\mathbf{z}, B\mathbf{w}\rangle=\mathbf{w}^{*}B^{*} J B\mathbf{z}.$$

Let $\mathbf{z}$ and $\mathbf{w}$ vary over a basis for
$\mathbf{C}^{2,1}$, we see that $B^{-1}=JB^{*}J.$ This means that
the inverse of $B\in\mathbf{PU}(2,1)$ has the following form:
$$B= \begin{bmatrix} a& b& c\\
d& e& f\\g&h& j\end{bmatrix}, \,\ B^{-1}= \begin{bmatrix} \overline{j}&\overline{ f}& \overline{c}\\
\overline{h}& \overline{e}& \overline{b}\\\overline{g}&\overline{d}&
\overline{a}\end{bmatrix}. \leqno(1)$$

We define the Siegel domain model of the complex hyperbolic 2-space,
$\mathbf{H}_{\mathbf{C}}^2$ as follows. We identify points of
$\mathbf{H}_{\mathbf{C}}^2$ with their horospherical coordinatess,
$z=(\xi,\nu,\mu)\in \mathbf{C}\times\mathbf{ R} \times
\mathbf{R}_{+}=\mathbf{H}_{\mathbf{C}}^2.$ Similarly, points in
$\partial\mathbf{H}_{\mathbf{C}}^2=\mathbf{C}\times \mathbf{R
}\times \{\infty\}$ are either $z=(\xi,\nu,0)\in
\mathbf{C}\times\mathbf{ R} \times \{0\}$ or a point at infinity,
denoted $q_{\infty}$. Define the map
$\psi:\overline{\mathbf{H}_{\mathbf{C}}^2}\rightarrow
\mathbf{PC}^{2,1}$ by
 $$\psi: (\xi,\nu,\mu)\mapsto \begin{bmatrix} -|\xi|^{2}-\mu+i\nu\\ \sqrt{2}\xi\\ 1\end{bmatrix} \quad  {\rm for}  \quad (\xi,\nu,\mu)\in \overline{\mathbf{H}_{\mathbf{C}}^2}-{q_{\infty}}, $$
and
 $$\psi:q_{\infty}\mapsto \begin{bmatrix} 1\\ 0\\ 0\end{bmatrix}.$$

 The map $\psi$ is a homeomorphism from $\mathbf{H}_{\mathbf{C}}^2$
 to the set of points $\mathbf{z}$ in $\mathbf{PC}^{2,1}$ with
 $\langle\mathbf{z},\mathbf{z}\rangle<0$. Also $\psi$ is a homeomorphism from
 $\partial\mathbf{H}_{\mathbf{C}}^2$ to  the set of points
 $\mathbf{z}$ with $\langle\mathbf{z},\mathbf{z}\rangle=0$. Let $L$ be a complex
 line intersecting $\mathbf{H}_{\mathbf{C}}^2$. Then $\psi(L)$ is a
 two-dimensional complex linear subspace of $\mathbf{C}^{2,1}$. The orthogonal
 complement of this subspace is a one (complex)-dimensional subspace
 of $\mathbf{C}^{2,1}$ spanned by a vector $\mathbf{p}$ with
 $\langle \mathbf{p}, \mathbf{p} \rangle>0$. Without loss of generality, we take
 $\langle \mathbf{p}, \mathbf{p}\rangle=1$ and call $\mathbf{p}$ the {\it polar
 vector}
 corresponding to the complex line $L$ (see page 75 of \cite{go}). The
 {\it Bergman metric} on $\mathbf{H}_{\mathbf{C}}^2$ is defined by the following
 formula for distance $\rho$ between points $z$ and
 $w$ of $\mathbf{C}^{2,1}$

      $$\cosh\big(\rho(z,w)/2\big)=\frac{\langle\psi(z),\psi(w)\rangle \langle\psi(w),\psi(z)\rangle}{\langle\psi(z),\psi(z)\rangle \langle\psi(w),\psi(w)\rangle}.$$

As in real hyperbolic geometry, $A$ holomorphic complex hyperbolic
isometry $g$ is said to be:

 (i)  {\it loxodromic} if it fixes no point  in $\mathbf{H}_{\mathbf{C}}^2$ but exactly two points of $\partial
\mathbf{H}_{\mathbf{C}}^2$;

(ii)  {\it parabolic} if it fixes fixes no point  in
$\mathbf{H}_{\mathbf{C}}^2$ but exactly one point of $\partial
 \mathbf{H}_{\mathbf{C}}^2$;

(iii)  {\it elliptic} if it fixes at least one point of $
 \mathbf{H}_{\mathbf{C}}^2$.

  The matrices corresponding to a loxodromic element and a parabolic element can be found in [10]. We
   will only give some matrices corresponding to the elliptic elements in this paper.
If $A $ is a elliptic element, then there are now three cases.
First, suppose that $A$ has a repeated eigenvalue with a two
dimensional eigenspace containing both positive and negative
vectors. This eigenspace corresponds to a complex line $L$ on which
$A$ acts as the identity. In particular, there are points of $
 \partial\mathbf{H}_{\mathbf{C}}^2$
fixed by $A$ and so $A$ is called boundary elliptic. As $A$ fixes
$L$ and rotates $
 \mathbf{H}_{\mathbf{C}}^2$ around $L$, it is complex reflection in the line
$L$. If $A$ is not boundary elliptic, then it has an eigenspace
spanned by a negative vector $\mathbf{w}$. This corresponds to a
fixed point $w\in \mathbf{H}_{\mathbf{C}}^2$. In this case $A$ is
called regular elliptic. There are two possibilities. Either $A$ has
a repeated eigenvalue with an eigenspace spanned by two positive
vectors. In this case $A$ is complex reflection in the point $w$.
Otherwise, $A$ has three distinct eigenvalues.

\begin{pro}  {\rm (1)} {\it Every boundary elliptic element  in $\mathbf{H}_{\mathbf{C}}^2$ is conjugate to
$$  \begin{bmatrix} u^{-{1/3}}& 0 & 0 \cr 0
&u^{2/3} & 0 \cr 0& 0
    &u^{-{1/3}}\end{bmatrix},
$$
    where} $u=e^{i\theta}$.

  {\rm (2)} {\it Every regular elliptic element  in
$\mathbf{H}_{\mathbf{C}}^2$ is conjugate to
$$  \begin{bmatrix} (u+w)/2& 0 & (u-w)/2 \cr 0
&v & 0 \cr (u-w)/2& 0
    &(u+w)/2\end{bmatrix},
$$
    where} $w=e^{2i\theta/3}$ and $u=v=e^{-i\theta/3}$.
\end{pro}

 Suppose that $A\in \mathbf{SU}(2,1)$ is an  elliptic
element. We define the {\it order} of $A$ as
   $${\rm order}(A)=\inf\{m>0, A^{m}=I\}.$$

As in the case of real hyperbolic geometry, a discrete subgroup of
$\mathbf{SU}(2,1)$ can not contain elliptic elements of infinite
order.

\section{The heisenberg group}
 Just as the boundary of real hyperbolic space may be identified
 with the one point compactification  of Euclidean space, so the
 boundary of complex hyperbolic space may be identified with one
 point compactification of the Heisenberg group. We now collect some
 of the basic facts about the Heisenberg group that will be used
 later.

 Consider the  $3$ dimensional Heisenberg group $\mathfrak{R}$ which is the set $\mathbf{C}\times
 \mathbf{R}$ (with coordinatess $(\xi,\nu))$ endowed with the
 multiplication law
  $$(\xi_{1},\nu_{1})\diamond
  (\xi_{2},\nu_{2})=(\xi_{1}+\xi_{2},\nu_{1}+\nu_{2}+2\Im\langle\langle\xi_{1},\xi_{2}\rangle\rangle),$$
where $\langle\langle\cdot,\cdot\rangle\rangle$ is the standard
positive definite
  Hermitian form on $\mathbf{C}$.
  The Heisenberg norm assigns to $(\xi,\nu)$ the non-negative real
  number

  $$|(\xi,\nu)|_{0}=(\|\xi\|^{4}+\nu^{2})^{\frac{1}{4}}=|\|\xi\|^{2}-i\nu|^{\frac{1}{2}}$$
where
$\|\xi\|^{2}=\langle\langle\xi,\xi\rangle\rangle=\sum|\xi_{i}|^{2}$.
This enables us to define the {\it Cygan metric} on the Heisenberg
group:
$$\rho_{0}((\xi_{1},\nu_{1}),(\xi_{2},\nu_{2}))=|(\xi_{1}-\xi_{2},\nu_{1}-\nu_{2}+2\Im\langle\langle\xi_{1},\xi_{2}\rangle\rangle)|_{0}=|(\xi_{1},\nu_{1})^{-1}\diamond
  (\xi_{2},\nu_{2})|_{0}.$$

  The Heisenberg group acts on itself by  Heisenberg translation.
  For $(\xi_{0},\nu_{0})\in \mathfrak{R}$, this is
$$T_{\xi_{0},\nu_{0}}: (\xi,\nu)\longmapsto
(\xi+\xi_{0},\nu+\nu_{0}+2\Im\langle\langle\xi_{0},\xi\rangle\rangle)=(\xi_{0},\nu_{0})\diamond
  (\xi,\nu).$$

Heisenberg group translation by $(0^{'},\nu_{0})$ where $0^{'}$ is
origin in $\mathbf{C}$ and $\nu_{0}\in \mathbf{R} $ are called
{\it vertical translations}.

\section{The ford isometric spheres}

In \cite{go} Goldman had extended the definition of isometric
spheres of M\"{o}bius transformations acting on the upper half space
to the Ford isometric spheres of complex hyperbolic transformations
of  the Siegel domain. These spheres and their associated geometric
properties have been extensively used in \cite{go,jp,p1,p2}.

Let $q_{\infty}=(1, 0, 0) \in \mathbf{C}^{2,1}$.

\begin{defn}\cite{p1} Let $X\in \mathbf{PU}(2,1)$. Suppose that $X$
does not fix $q_{\infty}$. Then the isometric sphere of $X$ is the
hypersurface

$$I_{X}=\{z\in \mathbf{H}_{\mathbf{C}}^2:
|\langle Z,q_{\infty}\rangle|=|\langle Z,X^{-1}(q_{\infty})\rangle |\}.$$ for any $Z\in
\mathbf{C}^{3}$ which maps onto $z$ projectively.
\end{defn}

 As in real case, $X$ maps $I_{X}$ to  $I_{X^{-1}}$ and  $X$ maps the component
of $\overline{\mathbf{H}_{\mathbf{C}}^2}\setminus I_{X}$ containing
$q_{\infty}$ to the component of
$\overline{\mathbf{H}_{\mathbf{C}}^2}\setminus I_{X^{-1}}$ not
containing $q_{\infty}$.
\begin{pro} \cite{p1}
 Let $X\in \mathbf{PU}(2,1)$ such that $X(q_{\infty})\neq q_{\infty}$. Then  the isometric
 sphere is the sphere for Cygan metric $\rho_{0}$ with center at
$X^{-1}(q_{\infty})$ and radius $r_{X}=\sqrt{\frac{1}{|g|}}.$
\end{pro}

\section{Proof of the theorem 3}
 In this section, we prove  Theorem 3. The basic structure of this proof resembles
the  original proof of \cite{gmm}. On the first hand,  we give a
proof of part one of Theorem 3.

Suppose that boundary elliptic element $A\in \mathbf{PU}(2,1)$ fixes
$0$ and $\infty$. This meas that complex line  $L_{A}$ fixed by $A$
is  spanned by $0$ and $\infty$. In other words
$$p_{A}=\begin{bmatrix} 0\\ 1\\ 0\end{bmatrix}.$$

Assume that $m\geq n$, let $f$ and $g$ be boundary elliptic elements
in $\mathbf{PU}(2,1)$  of order $m$ and $n$ respectively, that is,
$f$ and $g$ are reflections in complex lines and set
$$\delta=\delta(f,g)$$
and $\omega^{2}=e^{\delta+i\phi}$, where $\delta$ is the distance
between the complex lines fixed by $f$ and $g$ respectively.

The statement is invariant with respect to conjugation by elements
in $\mathbf{PU}(2,1)$. Thus by means of conjugation we may choose
some matrices form of $f$ and $g$ for the convenience of our
calculations.

We begin with the following two elements in $\mathbf{SU}(2,1)$

$$U_{1}= \begin{bmatrix} u_{1}^{-\frac{1}{3}}& 0& 0\\
0& u_{1}^{\frac{2}{3}}& 0\\0&0& u_{1}^{-\frac{1}{3}}\end{bmatrix}$$
and
$$U_{2}= \begin{bmatrix} u_{2}^{-\frac{1}{3}}& 0& 0\\
0& u_{2}^{\frac{2}{3}}& 0\\0&0& u_{2}^{-\frac{1}{3}}\end{bmatrix}$$
where $u_{1}=e^{\frac{2\pi i}{m}}, u_{2}=e^{\frac{2\pi i}{n}}.$

Obviously, $U_{1}$ and $U_{2}$ are boundary elliptic elements of
order $m$ and $n$. Moreover,  $U_{1}$ and $U_{2}$ fix the same
complex line with polar vector
$$p=\begin{bmatrix} 0\\ 1\\ 0\end{bmatrix}.$$

Now suppose that $A$ and  $B$ in $\mathbf{SU}(2,1)$ have the
following forms
$$A= \begin{bmatrix} \overline{\sqrt{\omega}}/2& \sqrt{\omega/2}& \overline{\sqrt{\omega}}2\\
1/\sqrt{2}& 0&
1/\sqrt{2}\\1/2\sqrt{\omega}&1/\overline{\sqrt{2\omega}}&-1/2\sqrt{\omega}
\end{bmatrix}$$ and
$$B= \begin{bmatrix} -1/2\overline{\sqrt{\omega}}& 1/\sqrt{2\omega}& 1/2\overline{\sqrt{\omega}}\\
1/\sqrt{2}& 0& 1/\sqrt{2}\\
\sqrt{\omega}/2& \overline{\sqrt{\omega/2}}&-\sqrt{\omega}/2
\end{bmatrix}.$$

We define the matrices representative $F$ and $G$ of $f$ and $g$ as
follows

$$F=AU_{1}A^{-1},\,\ G=BU_{2}B^{-1}.$$

Elementary calculations show that

$$F= \begin{bmatrix} \frac{1}{2}(u_{1}^{-\frac{1}{3}}+u_{1}^{\frac{2}{3}})& 0& \frac{|\omega|}{2}(u_{1}^{\frac{2}{3}}-u_{1}^{-\frac{1}{3}})\\
0& u_{1}^{-\frac{1}{3}}& 0\\
\frac{1}{2|\omega|}(u_{1}^{\frac{2}{3}}-u_{1}^{-\frac{1}{3}}) &0&
\frac{1}{2}(u_{1}^{-\frac{1}{3}}+u_{1}^{\frac{2}{3}})\end{bmatrix}$$
and
$$G= \begin{bmatrix} \frac{1}{2}(u_{2}^{-\frac{1}{3}}+u_{2}^{\frac{2}{3}})& 0& \frac{1}{2|\omega|}(u_{2}^{\frac{2}{3}}-u_{2}^{-\frac{1}{3}})\\
0& u_{2}^{-\frac{1}{3}}& 0\\
\frac{|\omega|}{2}(u_{2}^{\frac{2}{3}}-u_{2}^{-\frac{1}{3}}) &0&
\frac{1}{2}(u_{2}^{-\frac{1}{3}}+u_{2}^{\frac{2}{3}})\end{bmatrix}.$$

We can see that $f$ and $g$ are boundary elliptic elements of order
$m$ and $n$. The complex line fixed by $f$ has the polar vector
$p_{f}=(\sqrt{\omega/2}, 0, 1/\overline{\sqrt{2\omega}})^{T}$ and
the complex line fixed by $g$ has polar vector
$p_{g}=(\overline{\sqrt{\omega/2}}, 0, 1/{\sqrt{2\omega}})^{T}$.

\begin{pro}
Let $f$ and $g$ be boundary elliptic elements  in $\mathbf{PU}(2,1)$
 having the matrices of above. Then the complex lines fixed by $f$ and
$g$ with polar vectors $p_{f}$ and $p_{g}$ has distance $\delta$.
\end{pro}
\begin{pf}
Let the complex line $C_{f}$ fixed by $f$ with polar vector $p_{f}$
and the complex line $C_{g}$ fixed by $g$ with polar vector $p_{g}$.
Then by distance formulas in \cite{s}  we have

\begin{eqnarray*}
\begin {array}{ll}
{\rm dist}(C_{f},C_{g})&=2\cosh^{-1}(|\langle p_{f},p_{g}
\rangle|)\\&=2\cosh^{-1}(|\frac{\sqrt{\omega}}{\sqrt{2}}*\frac{\sqrt{\omega}}{\sqrt{2}}+
\frac{1}{\overline{\sqrt{2\omega}}}*\frac{1}{\overline{\sqrt{2\omega}}}|)\\&=
2\cosh^{-1}(|\frac{1}{2}(\omega +\frac{1}{\overline{\omega}})|)\\&=
\delta.
\end{array}
\end{eqnarray*}
\end{pf}

Suppose that $A\in \mathbf{SU}(2,1)$ does not fix $q_{\infty}$,
which  is equivalent to requiring that $g$ be non-zero when $A$ has
the form (1). Then the isometric sphere of $A$ is the sphere in the
Cygan metric with center $A^{-1}(\infty)$ and radius
$r_{A}=\frac{1}{\sqrt{|g|}}$. In Heisenberg coordinates

$$A^{-1}(\infty)=(\frac{\overline{h}}{\sqrt{2}\overline{g}},-
\Im\frac{j}{g}).$$ Similarly the  isometric sphere of $A^{-1}$ is
the Cygan sphere of radius $\frac{1}{\sqrt{|g|}}$ with center
$$A(\infty)=(\frac{d}{\sqrt{2}g},
\Im\frac{a}{g}).$$

The isometric spheres of $f$ and $g$ are easily calculated from
their matrices representative $F$ and $G$.  Using $u_{1}=e^{2\pi
i/m}$ we  have $u_{1}^{2/3}- u_{1}^{-1/3}=2ie^{\pi i/6m}\sin(\pi/m)$
and $u_{1}^{2/3}+u_{1}^{-1/3}=2e^{\pi i/6m}\cos(\pi/m)$. The Cygan
isometric sphere $I_{f}$ of $f$ has radius
$$r_{f}=\frac{1}{|\frac{1}{2|\omega|}(u_{1}^{\frac{2}{3}}-u_{1}^{-\frac{1}{3}})|}=\sqrt{\frac{|\omega|}{\sin(\pi/m)}}$$
and center
$$\bigg(0,
-\Im\frac{\frac{1}{2}(u_{1}^{-\frac{1}{3}}+u_{1}^{\frac{2}{3}})}{\frac{1}{2|\omega|}(u_{1}^{\frac{2}{3}}
-u_{1}^{-\frac{1}{3}})}\bigg)=\bigg(0,\frac{|\omega|\cos(\pi/m)}{\sin(\pi/m)}\bigg).$$
Similarly the Cygan isometric sphere $I_{f^{-1}}$ of $f^{-1}$ has
radius $r_{f^{-1}}=r_{f}$ and  center
$$\bigg(0,-\frac{|\omega|\cos(\pi/m)}{\sin(\pi/m)}\bigg).$$

 The Cygan isometric spheres $I_{g}$  and $I_{g^{-1}}$ of $g$
 have
radius
$$ r_{g}=r_{g^{-1}}=\frac{1}{\sqrt{\frac{|\omega|}{2}(u_{2}^{\frac{2}{3}}-u_{2}^{-\frac{1}{3}})}}=\sqrt{\frac{1}{|\omega|\sin(\pi/n)}}. $$
and the centers of  $I_{g}$  and $I_{g^{-1}}$  are

$$\bigg(0, \frac{\cos(\pi/n)}{|\omega|\sin(\pi/n)}\bigg),
\quad \bigg(0, -\frac{\cos(\pi/n)}{|\omega|\sin(\pi/n)}\bigg)$$
respectively.

The fundamental domain for the action of $f$ on the Heisenberg group
is the exterior of these two spheres $I_{f}$  and $I_{f^{-1}}$
together with the region bounded by their intersection.

We observe that the fundamental domain  of $f$ contains the
Heisenberg sphere $S_{f}^{*}$
 with center $(0,0)$ and radius $$ r_{f}^{\ast}=\sqrt{\frac{|\omega|}
{\sin(\pi/m)}\big(1-\cos(\pi/m)\big)}.$$
 The  Cygan isometric
spheres $I_{g}$  and $I_{g^{-1}}$ of $g$ are contained  in the
Heisenberg sphere $S_{g}^{*}$ with center $(0,0)$ and radius $$
r_{g}^{\ast}=\sqrt{\frac{1}
{|\omega|\sin(\pi/n)}\big(1+\cos(\pi/n)\big)}.$$

The interiors of $I_{g}$  and $I_{g^{-1}}$ are contained in the
interiors of $I_{f}$  and $I_{f^{-1}}$ if $r_g^{*}\leq r_f^{*}$.
That is
$$|\omega|^{2}=e^{\delta}
\geq\frac{\sin(\pi/m)}{1-\cos(\pi/m)}\frac{1+\cos(\pi/n)}{\sin(\pi/n)}.$$

Using
$$\frac{\sin\theta}{1-\cos\theta}=\frac{1+\cos\theta}{\sin\theta}$$

This translates into
$$\cosh(\delta)=\frac{|\omega|^{2}+|\omega|^{-2}}{2}\geq \frac{\cos\frac{\pi}{m}\cos\frac{\pi}{n}+1}{\sin\frac{\pi}{m}\sin\frac{\pi}{n}}.$$

 We have therefore seen that the exterior of a fundamental domain
for $\langle g \rangle$ lies inside a fundamental domain for
$\langle f \rangle$. It follows from the simplest version of
Kleinian-Maskit combination theorem that the group $\langle f, g
\rangle$ is discrete and isomorphic to the free product of cyclic
groups, $$\langle f,g\rangle\cong \langle f\rangle \ast \langle
g\rangle \cong \mathbf{Z}_{m}*\mathbf{Z}_{n}. $$

It is straightforward to extend the main result to the case where
either or both of $f$ and $g$ are complex reflection in a point. If
$f$ and $g$ are complex reflections in a point then the expressions
for $F$ and $G$ on the above become

$$F= \begin{bmatrix} \frac{1}{2}(u_{1}^{-\frac{1}{3}}+u_{1}^{\frac{2}{3}})& 0& \frac{|\omega|}{2}(u_{1}^{-\frac{1}{3}}-u_{1}^{-\frac{2}{3}})\\
0& u_{1}^{-\frac{1}{3}}& 0\\
\frac{1}{2|\omega|}(u_{1}^{-\frac{1}{3}}-u_{1}^{\frac{2}{3}}) &0&
\frac{1}{2}(u_{1}^{-\frac{1}{3}}+u_{1}^{\frac{2}{3}})\end{bmatrix}$$
which fixes $p_{f}=(\sqrt{\omega/2}, 0,
1/\overline{\sqrt{2\omega}})^{T}$  and

$$G= \begin{bmatrix} \frac{1}{2}(u_{2}^{-\frac{1}{3}}+u_{2}^{\frac{2}{3}})& 0& \frac{1}{2|\omega|}(u_{2}^{\frac{2}{3}}-u_{2}^{-\frac{1}{3}})\\
0& u_{2}^{-\frac{1}{3}}& 0\\
\frac{|\omega|}{2}(u_{2}^{\frac{2}{3}}-u_{2}^{-\frac{1}{3}}) &0&
\frac{1}{2}(u_{2}^{-\frac{1}{3}}+u_{2}^{\frac{2}{3}})\end{bmatrix}.$$
which fixes $p_{g}=(- 1/\sqrt{2\omega}, 0,
\overline{\sqrt{\omega/2}})^{T}.$

The distance between the fixed points or lines may be calculated as
in \cite{s}.  Namely, when one of  $p_{f}$ and $p_{g}$ is in $V_{+}$
and the other in $V_{-}$( that is one of $f$ and $g$ is complex
reflection in a point and the other is complex  reflection in a
complex line) then the distance between this point and complex line
is $\delta(f,g)$ where

  $$\sinh^{2}\big(\frac{\delta(f,g)}{2}\big)=\frac{\langle p_{f}, p_{g}\rangle \langle p_{g}, p_{f}\rangle}{-\langle p_{f}, p_{f}\rangle \langle p_{g},
  p_{g}\rangle}=|\omega/2-1/2\overline{\omega} |^{2}.$$

Similarly, when  $p_{f}$ and $p_{g}$ are both in $V_{-}$, so $f$ and
$g$ each are complex reflection in a point then the distance between
these points is $\delta(f,g)$ where

 $$\cosh^{2}\big(\frac{\delta(f,g)}{2}\big)=\frac{\langle p_{f}, p_{g}\rangle \langle p_{g}, p_{f}\rangle}{\langle p_{f}, p_{f}\rangle \langle p_{g},
  p_{g}\rangle}=|\omega/2+1/2\overline{\omega} |^{2}.$$

In either case

         $$\cosh^{2}\big(\frac{\delta(f,g)}{2}\big)=\frac{|\omega|^{2}+|\omega|^{-2}}{2}.$$

The same identity holds in the case where $f$ and $g$ fix complex
lines and $\delta(f,g)$ denotes the distance between these complex
lines.

 In each case the isometric spheres and fundamental domains are the
 same and so the other calculations go through with no changes.

\begin{rmk}
  The group generated by $f$ and $g$ preserves a
(unique) complex line $L$.  The restriction of the Bergman
metric to $L$ is just the Poincar\'{e} metric and both $f$ and $g$ act on $L$ as elliptic
hyperbolic isometries. Theorem 3 is a natural generalisation of
the result for real hyperbolic space of dimesions 2 and 3.
\end{rmk}

 \vskip 12pt {\bf Acknowledgement}. The  first author wishes to
thank Dr. Ying  Zhang and Dr.  Ser Peow Tan for their help and
support. This research was supported by National Natural Science
Foundational of China (No. 10671059) and the  first author also
supported by Chuang Xin Ji Jin of Hunan Province(No. 521298290).

\end{document}